\documentclass[11pt,a4paper]{amsart}
\usepackage[all]{xy}

\newcommand{\ie}{{\em i.e.}\ }
\newcommand{\cf}{{\em cf.}\ }

\newcommand{\ko}{\: , \;}
\newcommand{\ul}[1]{\underline{#1}}
\newcommand{\ol}[1]{\overline{#1}}

\newtheorem*{theorem}{Theorem}
\newtheorem*{lemma}{Lemma}
\newtheorem*{proposition}{Proposition}
\newtheorem*{corollary}{Corollary}

\newcommand{\opname}[1]{\operatorname{\mathsf{#1}}}

\renewcommand{\mod}{\opname{mod}\nolimits}

\newcommand{\Mod}{\opname{Mod}\nolimits}

\newcommand{\ind}{\opname{ind}\nolimits}
\newcommand{\per}{\opname{per}\nolimits}
\newcommand{\add}{\opname{add}\nolimits}
\newcommand{\op}{^{op}}

\newcommand{\colim}{\opname{colim}}
\newcommand{\cok}{\opname{cok}\nolimits}

\renewcommand{\ker}{\opname{ker}\nolimits}

\newcommand{\Z}{\mathbb{Z}}

\newcommand{\iso}{\stackrel{_\sim}{\rightarrow}}
\newcommand{\id}{\mathbf{1}}

\newcommand{\Hom}{\opname{Hom}}

\newcommand{\Ext}{\opname{Ext}}

\newcommand{\ca}{{\mathcal A}}

\newcommand{\cc}{{\mathcal C}}
\newcommand{\cd}{{\mathcal D}}
\newcommand{\ce}{{\mathcal E}}

\newcommand{\ch}{{\mathcal H}}
\newcommand{\ci}{{\mathcal I}}

\newcommand{\cm}{{\mathcal M}}

\newcommand{\cp}{{\mathcal P}}

\newcommand{\cs}{{\mathcal S}}
\newcommand{\ct}{{\mathcal T}}
\newcommand{\cu}{{\mathcal U}}
\newcommand{\cv}{{\mathcal V}}

\newcommand{\eps}{\varepsilon}

\newcommand{\cmp}{\opname{CMP}}
\newcommand{\cmi}{\opname{CMI}}

\newcommand{\cmac}{\opname{CM}}

\newcommand{\bt}{\circ}

\begin{document}

\title[Cluster tilted algebras]{Cluster-tilted
algebras are Gorenstein\\
and stably Calabi-Yau}

\author{Bernhard Keller}
\address{UFR de Math\'ematiques, UMR 7586 du CNRS, Case 7012,
   Universit\'e Paris 7, 2 place Jussieu, 75251 Paris Cedex 05,
   France}
\email{keller@math.jussieu.fr}
\author{Idun Reiten}
\address{Institutt for matematiske fag, Norges
Teknisk-naturvitenskapelige universitet, N-7491, Trondheim,
Norway} \email{idunr@math.ntnu.no}

 \date{December 20, 2005,
last modified on November 5, 2006}


\begin{abstract} We prove that in a $2$-Calabi-Yau triangulated
category, each cluster tilting subcategory is Gorenstein with all its
finitely generated projectives of injective dimension at most one.  We
show that the stable category of its Cohen-Macaulay modules is
$3$-Calabi-Yau. We deduce in particular that cluster-tilted algebras
are Gorenstein of dimension at most one, and hereditary if they are of
finite global dimension. Our results also apply to the stable (!)
endomorphism rings of maximal rigid modules of
\cite{GeissLeclercSchroeer05a}.  In addition, we prove a general
result about relative $3$-Calabi-Yau duality over non stable
endomorphism rings. This strengthens and generalizes the Ext-group
symmetries obtained in \cite{GeissLeclercSchroeer05a} for simple
modules. Finally, we generalize the results on relative Calabi-Yau
duality from $2$-Calabi-Yau to $d$-Calabi-Yau categories.  We show how
to produce many examples of $d$-cluster tilted algebras.
\end{abstract}

\maketitle


\section{Introduction}
\label{s:Introduction}

Let $k$ be a field and $H$ a finite-dimensional hereditary algebra.
The associated cluster category $\cc_H$ was introduced in
\cite{BuanMarshReinekeReitenTodorov04} and, for algebras $H$ of
Dynkin type $A_n$, in \cite{CalderoChapotonSchiffler04}.
It serves in the representation-theoretic approach to
cluster algebras introduced and studied by Fomin-Zelevinsky in
a series of articles \cite{FominZelevinsky02} \cite{FominZelevinsky03}
\cite{BerensteinFominZelevinsky05}. We refer to
\cite{FominZelevinsky03a} for more background on cluster algebras and
to
\cite{AssemBruestleSchiffler06}
\cite{AssemBruestleSchifflerTodorov05}
\cite{BuanMarshReiten05}
\cite{BuanMarshReitenTodorov05}
\cite{BuanReiten05a}
\cite{BuanReiten05b}
\cite{BuanReitenSeven06}
\cite{CalderoKeller05a}
\cite{CalderoKeller05b}
\cite{GeissLeclercSchroeer05c}
\cite{Iyama05}
\cite{IyamaReiten06}
\cite{Ringel06}
\cite{Thomas05}
for some
recent developments in the study of their links with representation theory.

The cluster category $\cc_H$ is the orbit category of the bounded
derived category of finite-dimensional right $H$-modules under the
action of the unique autoequivalence $F$ satisfying
\[
D\Hom(X, SY) = \Hom(Y,SFX) \ko
\]
where $D$ denotes the duality functor $\Hom_k(?,k)$ and $S$ the
suspension functor of the derived category. The functor $F$ is
$\tau S^{-1}$, where $\tau$ is the Auslander-Reiten translation
in the bounded derived category, or alternatively, $F$ is
$\Sigma S^{-2}$, where $\Sigma$ is the Serre
functor of the bounded derived category.
The cluster category is triangulated and Calabi-Yau
of CY-dimension $2$, \cf \cite{Keller05}. We still denote by $S$ its
suspension functor. A cluster-tilted algebra \cite{BuanMarshReiten04}
is the endomorphism algebra of a cluster tilting object $T$ of
$\cc_H$; this means that $T$ has no self-extensions, \ie
\[
\Hom(T,ST)=0 \ko
\]
and any direct sum $T\oplus T'$ with an indecomposable $T'$
not occuring as a direct summand of $T$ does have selfextensions.

We will prove that each cluster-tilted algebra is Gorenstein in the
sense of \cite{AuslanderReiten91}, \cf also \cite{Happel91},
and that its projective objects are of
injective dimension $\leq 1$ (hence its injective objects are of
projective dimension $\leq 1$). We deduce that if it has finite global
dimension, then it is hereditary. In particular, if it is given as the
quotient of a directed quiver algebra by an admissible ideal, then this ideal
has to vanish. We will also show that the stable category of
Cohen-Macaulay modules over a cluster-tilted algebra is
Calabi-Yau of CY-dimension $3$.

We work in a more general framework which also covers the stable (!)
endomorphism algebras of the maximal rigid modules of
\cite{GeissLeclercSchroeer05a}.
In our framework, the triangulated
category $\cc_H$ is replaced by an arbitrary triangulated category
$\cc$ with finite-dimensional $\Hom$-spaces and which is Calabi-Yau of
CY-dimension $2$, \ie we have
\[
D\Hom(X,Y)=\Hom(Y,S^2X) \ko X,Y\in\cc.
\]
This holds for cluster categories (including those coming from
$\Ext$-finite hereditary abelian categories with a tilting object, \cf
\cite{BuanMarshReinekeReitenTodorov04} \cite{Keller05}, and also those
coming from the $\Ext$-finite hereditary abelian $k$-categories in
\cite{ReitenVandenBergh02}, \cite{LenzingReiten05}, see
\cite{Keller05}), for stable module categories of preprojective
algebras of Dynkin graphs (\cf \cite[3.1, 1.2]{AuslanderReiten96},
\cite[8.5]{Keller05}) as investigated in
\cite{GeissLeclercSchroeer05a} and also for stable categories
$\ul{\cmac}(R)$ of Cohen-Macaulay modules over commutative complete
local Gorenstein isolated singularities of dimension~$3$, using
\cite{Auslander76}. The tilting object is replaced by a full
$k$-linear subcategory $\ct\subset\cc$ such that, among other
conditions, we have
\[
\Hom(T,ST')=0
\]
for all $T,T'\in \ct$ and where $\ct$ is maximal for this property. We then
consider the category $\mod \ct$ of finitely presented (right) modules
over $\ct$ and show that it is abelian and Gorenstein, and that the
homological dimensions of the projectives and the injectives are
bounded by $1$.  We also prove that $\mod\ct$ is equivalent to the
quotient of $\cc$ by an ideal, thus extending a result of
\cite{BuanMarshReiten04}, \cf also \cite{Zhu05}, and we generalize
Theorem~4.2 of \cite{BuanMarshReiten04} to give a
relationship between $\mod\ct$ and $\mod\ct'$ for `neighbouring'
subcategories $\ct$ and $\ct'$ of the type we consider, where
$\cc$ is a Calabi-Yau category of CY-dimension~2 whose cluster tilting
subcategories have `no loops'.
We give examples where $\ct$ has an infinite number of
isomorphism classes of indecomposables.

In section~\ref{s:Calabi-Yau}, we show that the stable category of
Cohen-Macaulay modules over $\ct$ is Calabi-Yau of CY-dimension
$3$. In the case where $\ct$ is invariant under the Serre functor, this
result was first proved in \cite{GeissKeller05}. It closely
resembles the Ext-group symmetries proved in Proposition~6.2 of
\cite{GeissLeclercSchroeer05a} for simple modules over non stable
endomorphism rings. We strengthen and generalize these Ext-group symmetries
in section~\ref{s:relative}, where we prove a
relative $3$-Calabi-Yau duality
over non stable endomorphism rings.

In section~\ref{s:relative-dp1}, we generalize
the results obtained on relative Calabi-Yau duality
to $d$-Calabi-Yau categories (\cf also \cite{Wraalsen05},
\cite{Thomas05}).
We also show how to produce many examples of $d$-cluster tilting
subcategories in $d$-cluster categories. For further examples,
we refer to \cite{IyamaReiten06}.

\section*{Acknowledgments}

This article has grown out of discussions which the authors had during the
meeting `Interactions between noncommutative algebra and algebraic
geometry' held at the BIRS in Banff in September 2005.  We are
grateful to the organizers of the meeting and in particular to Colin
Ingalls and James Zhang. We thank Osamu Iyama for pointing out
a mistake in a previous version of this article and for
allowing us to include his example~\ref{ss:Iyamas-example}.
We are grateful to Xiao-Wu Chen for pointing out missing conditions in
section~\ref{ss:reminder-Gorenstein}.

\section{The Gorenstein property}
\label{s:Gorenstein-property}

\subsection{The main result} \label{ss:main}
Let $k$ be a field and $\cc$ a triangulated $k$-linear category
with split idempotents and suspension functor $S$. We suppose that
all $\Hom$-spaces of $\cc$ are finite-dimensional and that $\cc$ admits
a Serre functor $\Sigma$, \cf \cite{ReitenVandenBergh02}.
We suppose that $\cc$ is Calabi-Yau of
CY-dimension $2$, \ie there is an isomorphism of triangle functors
\[
S^2 \iso \Sigma.
\]
Note that we do not exclude the possibility that there might already exist
an isomorphism $S^d\iso \Sigma$ for $d=0$ or $d=1$.
We fix such an isomorphism once and for all.
For $X,Y\in \cc$ and $n\in\Z$, we put
\[
\Ext^n(X,Y)=\cc\,(X,S^n Y).
\]
Assume that $\ct\subset\cc$ is a {\em cluster tilting
subcategory}. By this, we mean that $\ct$ is a maximal
$1$-orthogonal subcategory in the sense of Iyama~\cite{Iyama05}, \ie
\begin{itemize}
\item[a)] $\ct$ is a $k$-linear subcategory,
\item[b)] $\ct$ is contravariantly finite in $\cc$, \ie the
functor $\cc(?,X)|\ct$ is finitely generated for all $X\in \cc$,
\item[c)] we have $\Ext^1(T,T')=0$ for all $T,T'\in \ct$ and
\item[d)] if $X\in \cc$ satisfies $\Ext^1(T,X)=0$
for all $T\in \ct$, then $X$ belongs to $\ct$.
\end{itemize}
Note that condition b) means that for each $X\in \cc$, there is a
{\em right $\ct$-approximation}, \ie a morphism $T\to X$ with $T\in \ct$
such that each morphism $T' \to X$ with $T'\in \ct$ factors
through $T$. Condition d) is self-dual (because
of the Calabi-Yau property) and so are conditions a) and c). Part b)
of the proposition below shows that the dual of condition
b) also holds for $\ct$.

We point out that there are examples of Calabi-Yau categories of
CY-dimension~$2$ where we have cluster tilting subcategories with an
infinite number of isomorphism classes of indecomposables. They
arise from certain $\Ext$-finite hereditary abelian categories
with Serre duality from \cite{ReitenVandenBergh02} \cite{LenzingReiten05}.
For example, let
\[
\xymatrix{ \ldots \ar[r] & \circ \ar[r] & \circ\ar[r] & \circ\ar[r] & \ldots}
\]
be the $A^\infty_\infty$-quiver with linear orientation and $\ch$ the
category of its finite-dimensional representations over $k$. Let
$\cc=\cc_\ch$ be the corresponding cluster category. The
Auslander-Reiten quiver of $\ch$ is of the form $\Z A_\infty$.
Corresponding to a section of $\Z A_\infty$ with zigzag orientation,
there is a tilting class $\theta$ for $\ch$ in the sense of
\cite{LenzingReiten05}, that is
\begin{itemize}
\item[a)] $\Ext^1(T_1, T_2)=0$ for $T_1, T_2$ in $\theta$,
\item[b)] if, for $X\in \ch$, we have $\Hom(T,X)=0=\Ext^1(T,X)=0$ for
all $T$ in $\theta$, then $X=0$,
\item[c)] $\theta$ is locally bounded, \ie
for each indecomposable $T_1$ of $\theta$, there are only finitely many
indecomposables $T_2$ of $\theta$ such that $\Hom(T_1,T_2)\neq 0$ or
$\Hom(T_2,T_1)\neq 0$.
\end{itemize}
It is easy to see that $\theta$ induces a cluster tilting subcategory
$\ct$ in $\cc$.

By a {\em $\ct$-module}, we mean a contravariant $k$-linear functor
from $\ct$ to the category of vector spaces.
We denote by $\mod \ct$ the category of finitely presented
$\ct$-modules and by $F: \cc\to\mod\ct$ the functor which sends
$X\in\cc$ to the module $T \mapsto \Hom(T,X)$. Note that, since
idempotents split in $\ct$, this functor induces an equivalence
from $\ct$ to the category of projectives of $\mod \ct$. Following
\cite{AuslanderReiten91} and \cite{Happel91} we say that $\ct$
is {\em Gorenstein} if the finitely presented projective $\ct$-modules
are of finite injective dimension and the finitely presented
injective $\ct$-modules of finite projective dimension. Statements
b), b') and c) of the following theorem extend results of
\cite{BuanMarshReiten04}, \cf also \cite{Zhu05}.

\begin{proposition} \label{prop:gorenstein}
\begin{itemize}
\item[a)] The category $\mod \ct$ is abelian.
\item[b)] For each object $X\in \cc$, there is a triangle
\[
S^{-1} X \to T^X_1 \to T^X_0 \to X
\]
of $\cc$ with $T^X_i$ in $\ct$. In any such triangle, $T^X_0 \to X$ is
a right $\ct$-approximation and $X\to ST^X_1$ a left
$S\ct$-approximation.
\item[b')] For each object $X\in\cc$, there is a triangle
\[
X \to \Sigma T^0_X \to  \Sigma T^{1}_X \to S X
\]
in $\cc$ with $T^i_X$ in $\ct$. In any such triangle, $X \to \Sigma
T^0_X$ is a left $\Sigma\ct$-approximation and $\Sigma T^1_X \to SX$ a right
$\Sigma\ct$-approximation.
\item[c)] Denote by $(S\ct)$ the ideal of morphisms of $\cc$
which factor through an object of the form $ST$, $T\in\ct$.
Then the functor $F$ induces an equivalence $\cc/(S\ct)\iso \mod \ct$.
Moreover, in addition to inducing an equivalence from $\ct$ to the
category of projectives of $\mod\ct$, the functor $F$
induces an equivalence from $\Sigma\ct$ to the category of
injectives of $\mod \ct$.
\item[d)] The projectives of $\mod\ct$ are of injective dimension
at most $1$ and the injectives of projective dimension at most $1$.
Thus $\ct$ is Gorenstein of dimension at most $1$.
\end{itemize}
\end{proposition}

\begin{proof} a) Indeed, $\cc$ is a triangulated category. Thus
$\cc$ admits weak kernels. Now $\ct$ is contravariantly finite in
$\cc$ and therefore $\ct$ also admits weak kernels.
This is equivalent to the fact that $\mod \ct$ is abelian.

b) Since $\ct$ is contravariantly finite in $\cc$, there is a
right $\ct$-approximation $T_0 \to X$. Form a
triangle
\[
T_1  \to T_0 \to X \to S T_1 .
\]
For $T\in \ct$, the long exact sequence obtained by applying $\cc(T,?)$
to this triangle shows that $\Ext^1(T,T_1)$ vanishes. Thus $T_1$
lies in $\ct$. Conversely, if we are given such a triangle, then
$T_0 \to X$ is a right $\ct$-approximation because $\cc(T,ST_1)$
vanishes for $T\in \ct$, and $X \to ST_1$ is a left
$S\ct$-approximation because $\cc(T_0, ST)$ vanishes for
$T\in\ct$.

b') We apply (b) to $Y=S^{-1} X$ and obtain a triangle
\[
S^{-2}X \to T_1^Y \to T_0^Y \to S^{-1}X.
\]
Now we apply $S^2\iso \Sigma$ to this triangle to get the
required triangle
\[
X \to \Sigma T_1^Y \to \Sigma T_0^Y \to SX.
\]

c) Let $M$ be a finitely presented $\ct$-module. Choose a
presentation
\[
\cc(?,T_1) \to \cc(?,T_0) \to M \to 0.
\]
Form a triangle
\[
T_1 \to T_0 \to X \to ST_1.
\]
Since $\cc(T,ST_1)=0$, we obtain an exact sequence
\[
FT_1 \to FT_0 \to FX \to 0
\]
and $M$ is isomorphic to $FX$. One shows that $F$ is
full by lifting a morphism $FX \to FX'$ to a morphism between
presentations and then to a morphism of triangles.
Now let $X$, $X'$ be objects of $\cc$. Let $f: X \to X'$
be a morphism with $Ff=0$. By b), we have a triangle
\[
S^{-1}X \to T_1^X \to T_0^X \to X
\]
with $T_1^X$ and $T_0^X$ in $\ct$. Since
$\cc(T_0,f)=0$, the morphism $f$ factors through $ST_1$.

For $T\in\ct$, we have
\[
F\Sigma T = \cc(?, \Sigma T) = D\cc(T,?)
\]
so that $F\Sigma T$ is indeed injective. Since $\cc(ST', S^2 T)=0$ for
$T,T'\in\ct$, the functor $F$ induces a fully faithful functor
from $\Sigma\ct$ to the category of injectives of $\mod\ct$ and since
idempotents in $\cc$ split, this functor is an equivalence.

d) Given $X=\Sigma T$, $T\in\ct$, we form a triangle as
in b). We have
\[
FS^{-1} \Sigma T = F ST =0.
\]
This shows that the image under $F$ of the triangle of b) yields a
projective resolution of length $1$ of $F\Sigma T$.Similarly, to obtain an
injective resolution of $FT$, $T\in\ct$, we apply
b') to $X=T$ and use the fact that $FST=0$.
\end{proof}

\begin{corollary}
\begin{itemize}
\item[a)] Each $\ct$-module is either of infinite projective
(resp. injective) dimension or of projective (resp. injective)
dimension at most $1$.
\item[b)] The category $\mod\ct$ is either hereditary
or of infinite global dimension.
\end{itemize}
\end{corollary}

\begin{proof} a) Suppose that $M$ is a $\ct$-module of finite
projective dimension. Denote by $\Omega M$
the kernel of an epimorphism $P \to M$ with $P$ projective.
Suppose that $M$ is of projective dimension $n$. Then $\Omega^n M$
is projective. By induction on $p\geq 0$, we find that
$\Omega^{n-p}M$ is of injective dimension at most $1$.
In particular,
$M$ is of injective dimension at most $1$.
Dually, one shows the statement on the
injective dimensions. Part b) is immediate from a).
\end{proof}

\subsection{Comparing neighbours}
In this section, we briefly indicate how to generalize Theorem~4.2 of
\cite{BuanMarshReiten04} beyond cluster categories, namely to
Calabi-Yau categories of CY-dimension~$2$ whose cluster tilting
subcategories have `no loops'. This hypothesis holds not only for
cluster categories but also for the stable categories of
finite-dimensional preprojective
algebras~\cite{GeissLeclercSchroeer05a}.

We recall some results from \cite{Iyama05}. Let $\cc$ be a Calabi-Yau
category of CY-dimension~$2$. We assume that for each cluster tilting
subcategory $\ct$ of $\cc$ and for each indecomposable $T$ of $\ct$,
any non isomorphism $f: T \to T$ factors through an object $T'$ of
$\ct$ which does not contain $T$ as a direct summand. Now let
$\ct$ be a cluster tilting subcategory.  Then, according
to \cite{Iyama05}, the results on cluster categories from
\cite{BuanMarshReinekeReitenTodorov04} generalize literally:
For any indecomposable object $T$ of
$\ct$, there is an indecomposable $T^*$, unique up to isomorphism,
such that $T^*$ is not isomorphic to $T$ and the additive subcategory
$\ct'$ of $\cc$ with indecomposables
\[
\ind(\ct')=\ind(\ct) \setminus\{T\} \cup \{T^*\}
\]
is a cluster tilting subcategory of $\cc$. Moreover, there are
approximation triangles
\[
T^* \to B \to T \to ST^* \mbox{ and } T \to B' \to T^* \to ST
\]
connecting $T^*$ and $T$ such that $B$ and $B'$ belong to
both $\ct$ and $\ct'$. With these notations and assumptions, we get the
following connection between $\mod\ct$ and $\mod\ct'$.

\begin{proposition} Let $S_T$ and $S_{T^*}$ be the simple tops
of $\cc(?,T)$ and $\cc(?,T^*)$ respectively. Then there is an
equivalence of categories
\[
\mod\ct/ \add S_T \iso \mod \ct'/ \add S_{T^*}.
\]
\end{proposition}

The proof, left to the reader, follows the
proof of Theorem~4.2 of \cite{BuanMarshReiten04} and uses
the previous proposition.

\section{The Calabi-Yau property}
\label{s:Calabi-Yau}

\subsection{Reminder on exact categories} An {\em exact category} in the
sense of Quillen \cite{Quillen73} is an additive category $\ca$
endowed with a distinguished class of sequences
\[
\xymatrix{ 0 \ar[r] & A \ar[r]^i & B \ar[r]^p & C \ar[r] & 0} \ko
\]
where $i$ is a kernel of $p$ and $p$ a cokernel of $i$. We will
state the axioms these sequences have to satisfy using the
terminology of \cite{GabrielRoiter92}: The morphisms $p$ are called
deflations, the morphisms $i$ inflations and the pairs
$(i,p)$ conflations. The given class of conflations has to satisfy
the following natural axioms:
\begin{itemize}
\item[Ex0] The identity morphism of the zero object is a
deflation.
\item[Ex1] The composition of two deflations is a deflation.
\item[Ex2] Deflations are stable under base change, \ie if $p:
Y\to Z$ is a deflation and $f:Z'\to Z$ a morphism, then there is a
cartesian square
\[
\xymatrix{Y \ar[r]^p & Z \\
Y' \ar[u]^{f'} \ar[r]_{p'} & Z' \ar[u]_f}
\]
where $p'$ is a deflation.
\item[Ex2'] Inflations are stable under cobase change, \ie if
$i: X \to Y$ is an inflation and $g:X \to X'$ a morphism, there is
a cocartesian square
\[
\xymatrix{X \ar[d]^g \ar[r]^i & Y \ar[d]^{g'} \\
X' \ar[r]_{i'} & Y'}
\]
where $i'$ is an inflation.
\end{itemize}
As shown in \cite{Keller90}, these axioms are equivalent to Quillen's
and they imply that if $\ca$ is small, then there is a fully faithful
functor from $\ca$ into an abelian category $\ca'$ whose image is an
additive subcategory closed under extensions and such that a sequence
of $\ca$ is a conflation iff its image is a short exact sequence of $\ca'$.
Conversely, one easily checks that
an extension closed full additive subcategory $\ca$ of an
abelian category $\ca'$ endowed with all conflations which become
exact sequences in $\ca'$ is always exact.
The fundamental notions and constructions of homological algebra, and in
particular the construction of the derived category, naturally
extend from abelian to exact categories, \cf \cite{Neeman90} and
\cite{Keller96}.

A {\em Frobenius category} is an exact category $\ca$ which has enough
injectives and enough projectives and where the class of projectives
coincides with the class of injectives. In this case, the
{\em stable category} $\ul{\ca}$ obtained by dividing $\ca$ by
the ideal of morphisms factoring through a projective-injective
carries a canonical structure of triangulated category,
\cf \cite{Heller68} \cite{Happel87} \cite{KellerVossieck87}.
Its suspension functor $X \to SX$ is constructed by choosing,
for each object $X$ of $\ca$, a conflation
\[
X \to I \to SX \ko
\]
where $I$ is injective.

\subsection{Reminder on Gorenstein categories}
\label{ss:reminder-Gorenstein}
We refer to
\cite{AuslanderReiten91}
\cite{AuslanderReiten91a}
\cite{Happel91}
and the references given there for the
theory of Gorenstein algebras and their modules in non commutative and
commutative algebra. In this section, we restate some fundamental
results of the theory in our setup.

Let $\ca$ be a $k$-linear exact category with enough projectives
and enough injectives. 
Assume that $\ca$ is {\em Gorenstein}\footnote{We
thank Xiao-Wu Chen for pointing out that it is not enough to assume
that the projectives are of finite injective dimension and the injectives
of finite projective dimension.}, \ie
the full subcategory $\cp$ of the projectives is covariantly
finite, the full subcategory $\ci$ of injectives is contravariantly
finite and there is an integer $d$ such that all projectives are of injective
dimension at most $d$ and all injectives are of projective dimension
at most $d$. 
Let $\cd^b_f(\ca)$ denote the full triangulated subcategory of
the bounded derived category $\cd^b(\ca)$ generated by the
projectives (equivalently: the injectives).
Call an object $X$ of $\ca$ projectively (resp. injectively)
Cohen-Macaulay if it satisfies
\[
\Ext^i_\ca(X,P)=0 \mbox{ (resp. } \Ext^i_\ca(I,X)=0 \mbox{)}
\]
for all $i>0$ and all projectives $P$ (resp. injectives $I$).
Let $\cmp(\ca)$ resp. $\cmi(\ca)$ denote the full
subcategories of $\ca$ formed by the projectively resp.
injectively Cohen-Macaulay objects. Clearly, $\cmp(\ca)$ is an
exact subcategory of $\ca$. Endowed with this exact structure,
$\cmp(\ca)$ is a Frobenius category and its subcategory of
projective-injectives is $\cp$. Dually, $\cmi(\ca)$ is a Frobenius
category whose subcategory of projective-injectives is $\ci$.
We denote by $\ul{\cmp}(\ca)$ and $\ul{\cmi}(\ca)$ the stable
categories associated with these Frobenius categories.

In analogy with a classical result on Frobenius categories
(\cf \cite{KellerVossieck87}, \cite{Rickard89b}),
we have the following

\begin{theorem}[\cite{Happel91}] The inclusions
\[
\xymatrix{ \cmp(\ca) \ar[r] & \ca & \cmi(\ca) \ar[l] }
\]
induce triangle equivalences
\[
\xymatrix{ \ul{\cmp}(\ca) \ar[r] & \cd^b(\ca)/\cd^b_f(\ca) &
\ul{\cmi}(\ca) \ar[l]. }
\]
\end{theorem}
The {\em stable Cohen-Macaulay category} of $\ca$ is by definition
$\ul{\cmac}(\ca) = \cd^b(\ca)/\cd^b_{f}(\ca)$. It is a triangulated
category and an instance of a stable derived category in the sense of
\cite{Krause05}.

We can make the statement of the theorem more precise: The
canonical functor
\[
\ca/(\cp) \to \ul{\cmac}(\ca)
\]
admits a fully faithful left adjoint whose image is $\ul{\cmp}(\ca)$
and the canonical functor
\[
\ca/(\ci)\to \ul{\cmac}(\ca)
\]
admits a fully faithful right adjoint whose image is $\ul{\cmi}(\ca)$.

The stable Cohen-Macaulay category is also canonically equivalent to
the triangulated category $\Z \underline{\ca}$ obtained
\cite{Heller68} from the suspended \cite{KellerVossieck87} category
$\ul{\ca}=\ca/(\ci)$ by formally inverting its suspension functor $S$
and to the triangulated category obtained from the cosuspended
category $\ol{\ca}=\ca/(\cp)$ by formally inverting its loop functor
$\Omega$.  Thus, we have
\[
\ul{\cmac}(\ca)(X,S^n Y) = \colim_p \ul{\ca}(S^p X, S^{n+p} Y)
=\colim_p \ol{\ca}(\Omega^{p+n}X, \Omega^{p} Y)
\]
for all $X,Y\in\ca$.

We will need the following easily proved lemma (the statement of whose dual
is left to the reader).

\begin{lemma} \label{lemma:cm}
\begin{itemize}
\item[a)] Suppose that all injectives of $\ca$ are of projective
dimension at most $d$.
Then $S^d Y$ is injectively Cohen-Macaulay for all
$Y\in\ca$.
\item[b)] Let $n\geq 1$ and $X,Y\in\ca$. If $S^n Y$ is injectively
Cohen-Macaulay, we have
\[
\Ext^n_\ca(X,Y)/\Ext^n_\ca(I,Y) \iso \ol{\ca}(X,S^n Y)
\iso \ul{\cmac}(\ca)(X,S^n Y) \ko
\]
where $X \to I$ is an inflation into an injective. If
$S^{n-1}Y$ is injectively Cohen-Macaulay, then we have
\[
\Ext^n_\ca(X,Y) \iso \ol{\ca}(X,S^n Y) \iso \ul{\cmac}(\ca)(X, S^n Y).
\]
\end{itemize}
\end{lemma}

\subsection{Modules over cluster tilting subcategories} Let $\cc$
and $\ct$ be as in section~\ref{s:Gorenstein-property}. We have
seen there that $\mod \ct$ is a Gorenstein category in the sense
of the preceding section.

\begin{theorem} The stable Cohen-Macaulay category of $\mod\ct$ is
Calabi-Yau of CY-dimension~$3$.
\end{theorem}

Note that, by sections~\ref{ss:reminder-Gorenstein} and \ref{ss:main}, the theorem
implies that for all $\ct$-modules $X$, $Y$, we have a canonical
isomorphism
\[
D\Ext^2_{\mod\ct}(Y,X) \iso \ul{\Ext}^1_{\mod\ct}(X,Y) \ko
\]
where $\ul{\Ext}^1$ denotes the cokernel of the map
\[
\Ext^1_{\mod\ct}(I,Y) \to \Ext^1_{\mod\ct}(X,Y)
\]
induced by an arbitrary monomorphism $X \to I$ into an injective.

The theorem will be proved below in section~\ref{proof:cy3}.
If $\ct$ is stable under $\Sigma$, then $\mod\ct$ is a Frobenius
category. In this case, the theorem is proved in \cite{GeissKeller05}
using the fact that then $\ct$ is a `quadrangulated category' which is
Calabi-Yau of CY-dimension~$1$.

\subsection{Three simple examples} Consider the algebra $A$ given by
the quiver with relations
\[
\xymatrix@=0.5cm{ & & 3 \ar[dr]^\beta &  &\\
& 2 \ar[ru]^\alpha & & 4 \ar[ll]^\gamma &
\alpha\gamma=\beta\alpha=\gamma\beta=0. \\
1 \ar[ru]^\delta & & & &
}
\]
It is cluster-tilted of type $A_4$. For $1\leq i\leq 4$, denote by
$P_i$ resp. $I_i$ resp. $S_i$ the indecomposable projective, resp. injective
resp. simple (right) module associated with the vertex $i$. We have
$I_1=P_3$ and $I_3=P_4$ and we have minimal projective resolutions
\[
0 \to P_1 \to P_3 \to I_2 \to 0 \ko 0\to P_1 \to P_2 \to I_4 \to 0
\]
and minimal injective resolutions
\[
0 \to P_1 \to I_1 \to I_2 \to 0\ko 0 \to P_2 \to I_1\oplus I_4 \to I_2 \to 0.
\]
Thus the simple modules except $S_1$ are injectively Cohen-Macaulay
and the simple modules except $S_2$ are projectively Cohen-Macaulay.
Therefore, we have
\[
\Ext^1_A(S_i,S_j) \iso D\Ext^2_A(S_j,S_i)
\]
unless $(i,j)=(1,2)$. Now $\dim \Ext^1_A(S_i,S_j)$ equals
the number of arrows from $j$ to $i$ and $\Ext^2_A(S_j,S_i)$ is
isomorphic to the space of minimal relations from
$i$ to $j$, \cf \cite{Bongartz83}. Thus for example, the
arrow $\alpha$ corresponds to the relation $\gamma\beta=0$,
and similarly for the arrows $\beta$ and $\gamma$.
Note that no relation corresponds to the arrow $\delta$.
The stable category $\ul{\cmac}(\mod A)$ is equivalent
to the stable category of the full subquiver with relations
on $2,3,4$. It is indeed $3$-Calabi-Yau.

Consider the algebra $A$ given by the quiver with relations
\[
\xymatrix@=0.5cm{ & 2 \ar[dr]^\beta &  &\\
1 \ar[ru]^\alpha \ar[dr]_\delta & & 3 \ar[ll]^\eps &
\alpha\eps=\eps\beta=\delta\eps=\eps\gamma=0, \beta\alpha=\gamma\delta.\\
& 4 \ar[ru]_\gamma & & }
\]
It is cluster-tilted of type $D_4$. We have $P_1=I_3$ and $P_3=I_1$.
There are minimal injective resolutions
\[
0 \to P_2 \to I_1 \to I_4 \to 0 \ko
0 \to P_4 \to I_1 \to I_2 \to 0
\]
and minimal projective resolutions
\[
0 \to P_4 \to P_3 \to I_2 \to 0 \ko
0 \to P_2 \to P_3 \to I_4 \to 0.
\]
Thus the simples $S_1$ and $S_3$ are both projectively
and injectively Cohen-Macaulay and whenever $i,j$ are
connected by an arrow, we have
\[
\Ext^1_A(S_i,S_j) \iso D\Ext^2_A(S_j,S_i).
\]
Since there are neither arrows nor relations between
$2$ and $4$, we obtain a perfect correspondence between
arrows and relations: the outer arrows correspond to
zero relations and the inner arrow to the commutativity
relation.

Consider the algebra $A$ given by the quiver
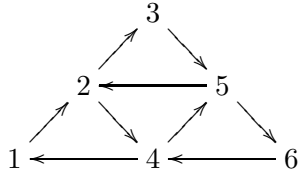
\[
\xymatrix@=0.5cm{ & & 3 \ar[dr] & & \\
 & 2 \ar[ru] \ar[rd] & & 5 \ar[dr] \ar[ll] & \\
1 \ar[ru] & & 4 \ar[ll] \ar[ru] & & 6 \ar[ll] }
\]
subject to the following relations: each path of length two
containing one of the outer arrows vanishes; each of the
three rhombi containing one of the inner arrows is commutative.
As shown in \cite{GeissLeclercSchroeer05b}, this algebra
is isomorphic to the stable endomorphism algebra of a
maximal rigid module over the preprojective algebra of
type $A_4$. It is selfinjective and thus Gorenstein.
According to the theorem, it is stably $3$-Calabi-Yau,
a fact which was first proved in \cite{GeissKeller05}
using the fact that the projectives over $A$ form a
quadrangulated category. Note that the isomorphisms
\[
\Ext^1_A(S_i,S_j) \iso D\Ext^2_A(S_j,S_i)
\]
translate into a perfect correspondence between arrows
and relations: each outer arrow corresponds to a zero
relation and each inner arrow to a commutativity relation.

\subsection{Proof of the theorem} \label{proof:cy3}
Denote by $\ul{\mod}\ct$ the
quotient of $\mod\ct$ by the ideal of morphisms factoring through
a projective and by $\ol{\mod}\ct$ the quotient by the ideal of
morphisms factoring through an injective. It follows from
Proposition~\ref{prop:gorenstein} that the functor $F$ induces
equivalences
\[
\cc/(\ct, S\ct) \iso \ul{\mod}\ct \mbox{ and }
\cc/(S\ct, \Sigma\ct) \iso \ol{\mod}\ct.
\]
Thus the suspension functor $S:\cc\to\cc$ induces a well-defined
equivalence $\tau : \ul{\mod}\ct \to \ol{\mod}\ct$. We first prove
that Serre duality for $\cc$ implies the Auslander-Reiten formula for
$\tau$ in $\mod\ct$. This shows in particular that $\tau$ is indeed
isomorphic to the Auslander-Reiten translation of $\mod \ct$. In the
context of cluster categories, it was proved in
\cite{BuanMarshReiten04} that the functor induced by $S$ coincides
with the Auslander-Reiten translation on objects of $\mod\ct$.
Note that if $\cc$ and $\mod\ct$ are Krull-Schmidt categories,
then, as for cluster-tilted algebras, we obtain the Auslander-Reiten
quiver of $\mod \ct$ from that of $\cc$ by removing the vertices
corresponding to indecomposables in $S\ct$.

\begin{lemma} Let $X,Y\in\mod\ct$. Then we have canonical
isomorphisms
\[
\ol{\Hom}(Y, \tau X) \iso D\Ext^1(X,Y) \iso
\ul{\Hom}(\tau^{-1}Y,X).
\]
\end{lemma}

\begin{proof} Choose $L$, $M$ in $\cc$ such that $FM=X$ and $FL=Y$.
Form triangles
\[
S^{-1}M \to T_1^M \to T_0^M \to M \mbox{ and } S^{-1}\Sigma T^1_L \to L
\to \Sigma T^0_L \to \Sigma T^1_L
\]
in $\cc$ as in Proposition~\ref{prop:gorenstein}. We obtain an exact sequence
\[
FS^{-1} M \to FT_1^M \to FT_0^M \to FM \to 0
\]
in $\mod\ct$ and its middle two terms are projective. Using this
sequence we see that $\Ext^1_\ct(FM,FL)$ is isomorphic to the
middle cohomology of the complex
\begin{equation} \label{eq:ext-one-seq}
\Hom(FT_0^M, FL) \to  \Hom(FT^M_1,FL) \to \Hom(FS^{-1}M,FL).
\end{equation}
Now if $V$ is an object of $\cc$ and
\[
V \to \Sigma T^0_V \to \Sigma T^1_V \to SV
\]
a triangle as in Proposition~\ref{prop:gorenstein}, then we have,
for each $U$ of $\cc$, an exact sequence
\[
\cc(U, ST^1_V) \to \cc (U, V) \to \Hom (FU, FV) \to 0.
\]
By applying this to the three terms in the
sequence~(\ref{eq:ext-one-seq}) and using the vanishing of $\Ext^1$
between objects of $\ct$, we obtain that $E=\Ext^1(FM, FL)$ is the
intersection of the images of the maps
\[
\cc(S^{-1}M, ST^1_L) \to \cc(S^{-1}M,L) \mbox{ and }
\cc(T_1^M,L) \to \cc(S^{-1}M,L).
\]
By the exact sequences
$\cc(T_1^M,L) \to \cc(S^{-1}M,L) \to \cc(S^{-1}T_0^M,L)$ and
\[
\cc(S^{-1}M, ST^1_L) \to \cc(S^{-1}M,L) \to \cc(S^{-1}M, \Sigma T^0_L)
\ko
\]
the group $E$ also appears as a kernel in the exact sequence
\[
0 \to E \to \cc(S^{-1}M, L ) \to \cc(S^{-1}M, \Sigma T^0_L) \oplus
\cc(S^{-1} T_0^M, L).
\]
If we dualize this sequence and use the isomorphism
\[
D\cc(U,V) \iso \cc(V,\Sigma U) =  \cc(V,S^2 U),
\]
we obtain the exact sequence
\[
\cc(\Sigma T^0_L, SM) \oplus \cc(L,S T_0^M) \to \cc(L,SM) \to DE
\to 0.
\]
This proves the left isomorphism. The right one follows since
$\tau$ is an equivalence.
\end{proof}

Now let $X,Y$ be $\ct$-modules and $\cs$ the stable
Cohen-Macaulay category of $\mod\ct$. We will construct a
canonical isomorphism
\[
D\cs(Y,S^2X)
\iso
\cs(X,SY).
\ko
\]
Since the injectives of $\mod\ct$ are of
projective dimension at most $1$ and the
projectives of injective dimension at most $1$, we have
\[
\Hom_{\cs}(Y,S^2X) = \Ext^2_{\mod\ct}(Y,X)
\]
and
\[
\Hom_{\cs}(X,SY) = \cok(\Ext^1_{\mod\ct}(I,Y) \to
\Ext^1_{\mod\ct}(X,Y))\ko
\]
by lemma \ref{lemma:cm} b),
where $X \to I$ is a monomorphism into an injective.
We will examine $\Ext^2_{\mod\ct}(Y,X)$ and compare its dual with
$\Ext^1_{\mod\ct}(X,Y)$ as described by the Auslander-Reiten formula.
Choose $M$ and $L$ in $\cc$ such that $X=FM$ and $Y=FL$. Form the triangle
\[
M \to \Sigma T^0_M \to \Sigma T^1_M \to SM.
\]
We obtain an exact sequence
\[
0 \to FM \to F\Sigma T^0_M \to F\Sigma T^1_M \to FSM \to FS\Sigma
T^0_M
\]
whose second and third term are injective. This is sufficient to
obtain that $\Ext^2(FL,FM)$ is isomorphic to the middle homology
of the complex
\[
\Hom(FL, F\Sigma T^1_M) \to \Hom(FL, FSM) \to \Hom(FL, FS\Sigma T^0_M).
\]
Now $F\Sigma T^1_M \to FSM$ is a right $\ci$-approximation and
$FM \to F\Sigma T^0_M$ is a monomorphism into an injective.
Put $I=F\Sigma T^0_M$. Then we have
\[
\Ext^2(FL, FM) = \ker (\ol{\Hom}(FL, \tau FM) \to \ol{\Hom}(FL,
\tau I)) \ko
\]
Using the Auslander-Reiten formula, we obtain
\[
D\Ext^2(FL,FM) \iso \cok(\Ext^1(I,FL) \to \Ext^1(FL,FM))
\]
and thus
\[
D\Ext^2(Y,X) \iso \cok(\Ext^1(I,Y) \to \Ext^1(X,Y)) \ko
\]
where $X\to I$ is a monomorphism into an injective. This
yields
\[
D\cs(Y,S^2 X) \iso \cs(X,SY)
\]
as claimed.

\section{Relative 3-Calabi-Yau duality over non stable endomorphism rings}
\label{s:relative}

Let $\ce$ be a $k$-linear Frobenius category with split idempotents.
Suppose that its stable category $\cc=\ul{\ce}$ has finite-dimensional
$\Hom$-spaces and is Calabi-Yau of CY-dimension $2$.
This situation occurs in the following examples:
\begin{itemize}
\item[(1)] $\ce$ is the category of finite-dimensional
modules over the preprojective algebra of a Dynkin quiver
as investigated in \cite{GeissLeclercSchroeer05a}.
\item[(2)] $\ce$ is the category of Cohen-Macaulay modules over a commutative
complete local Gorenstein isolated singularity of dimension~$3$.
\item[(3)] $\ce$ is a Frobenius category whose stable category
is triangle equivalent to the cluster category associated
with an $\Ext$-finite hereditary category. Such Frobenius
categories always exist by \cite{Keller05}, section~9.9.
\item[(4)] $\ce$ is a Frobenius category whose stable
category is triangle equivalent to the bounded derived category
of coherent sheaves on a Calabi-Yau surface (\ie a K3-surface).
For example, one obtains such a Frobenius category by taking
the full subcategory of the exact category of left bounded complexes
of injective quasi-coherent sheaves whose homology is bounded
and coherent.
\end{itemize}

Let $\ct\subset\cc$
be a cluster tilting subcategory and $\cm\subset\ce$ the preimage
of $\ct$ under the projection functor. In particular, $\cm$
contains the subcategory $\cp$ of the projective-injective
objects in $\cm$. Note that $\ct$ equals the quotient $\ul{\cm}$
of $\cm$ by the ideal of morphisms factoring through a projective-injective.

We know from section~\ref{s:Gorenstein-property} that $\mod \ul{\cm}$ is
abelian. In general, we cannot expect the category $\mod\cm$ of
finitely presented $\cm$-modules to be abelian (\ie $\cm$ to have weak
kernels). However, the category $\Mod\cm$ of all right $\cm$-modules
is of course abelian. Recall that the {\em perfect derived category}
$\per(\cm)$ is the full triangulated subcategory of the derived
category of $\Mod\cm$ generated by the finitely generated projective
$\cm$-modules. We identify $\Mod\ul{\cm}$ with the full subcategory of
$\Mod\cm$ formed by the modules vanishing on $\cp$.
The following proposition is based on the methods of \cite{Auslander66}.

\begin{proposition}
\begin{itemize}
\item[a)] For each $X\in\ce$, there is a conflation
\[
0 \to M_1 \to M_0 \to X \to 0
\]
with $M_i$ in $\cm$. In any such conflation, the morphism $M_0 \to X$
is a right $\cm$-approximation.
\item[b)] Let $Z$ be a finitely presented $\cm$-module. Then
$Z$ vanishes on $\cp$ iff there exists a conflation
\[
0 \to K \to M_1 \to M_0 \to 0
\]
with $M_0, M_1$ in $\cm$ such that $Z = \cok (M_1^\wedge \to M_0^\wedge)$,
where $M^\wedge$ denotes the $\cm$-module represented by $M\in\cm$.
In this case, we have
\[
Z \iso \Ext^1_\ce(?,K)|\cm.
\]
\item[c)] Let $Z$ be a finitely presented
$\ul{\cm}$-module. Then $Z$ considered as an $\cm$-module lies
in $\per(\cm)$ and we have a canonical isomorphism
\[
D\per(\cm)(Z,?) \iso \per(\cm)(?,Z[3]).
\]
\end{itemize}
\end{proposition}

\begin{proof}
a) Indeed, we know that there is a triangle
\[
M_1 \to M_0 \to X \to SM_1
\]
in $\ul{\ce}$ with $M_1$, $M_0$ in $\ct$. We lift it to the
required conflation.

b) Clearly, if $Z$ is of the form given, it vanishes on $\cp$.
Conversely, suppose that $Z$ is finitely presented with
presentation
\[
\xymatrix{M_1^\wedge \ar[r]^{p^\wedge} &  M_0^\wedge \ar[r] &  Z \ar[r] &  0}
\]
and vanishes on $\cp$. Then each morphism $P \to M_0$ with $P\in\cp$
lifts along $p$. Now let $q: P_0 \to M_0$ be a deflation with
projective $P_0$. Then we have $q=pq'$ for some $q'$.
Clearly, the morphism
\[
M_1 \oplus P_0 \to M_0
\]
with components $p$ and $q$ is a deflation and $M_1 \to M_0$ is
a retract of this deflation: The retraction is given by $[\id, q']$
and the identity of $M_0$. Since $\ce$ has split idempotents,
it follows that $M_1 \to M_0$ is a deflation. The last assertion
follows from the vanishing of $\Ext^1_\ce(M,M_1)$ for all $M$ in $\cm$.

c) Clearly, each representable $\ul{\cm}$-module is finitely presented as
an $\cm$-module. Since a cokernel of a morphism between finitely
presented modules is finitely presented, the module $Z$ considered
as an $\cm$-module is finitely presented. So we can choose a conflation
\[
0 \to K \to M_1 \to M_0 \to 0
\]
as in b). Now choose a conflation
\[
0 \to M_3 \to M_2 \to K \to 0
\]
with $M_i$ in $\cm$ as in a). Then the image of the spliced complex
\[
0 \to M_3 \to M_2 \to M_1 \to M_0 \to 0
\]
under the Yoneda functor $\cm \to \Mod \cm$ is a projective
resolution  $P$ of $Z$. So $Z$ belongs to $\per(\cm)$. We now prove the
duality formula.
We will exhibit a canonical linear form $\phi$ on $\per(\cm)(P,P[3])$
and show that it yields the required isomorphism. We have
\[
\per(\cm)(P,P[3])= \per(\cm) (P,Z[3]) =
\cok(\Hom(M_2^\wedge, Z) \to \Hom(M_3^\wedge, Z))
\]
and this is also isomorphic to
\[
\cok(Z(M_2) \to Z(M_3)) = \cok(\Ext^1_\ce(M_2,K) \to \Ext^1_\ce(M_3,K)).
\]
Now the conflation
\[
0 \to M_3 \to M_2 \to K \to 0
\]
yields an exact sequence
\[
\Ext^1_\ce(M_2, K) \to \Ext^1_\ce(M_3,K) \to
\Ext^2_\ce(K,K).
\]
Thus we obtain an injection
\[
\per(\cm)(P,P[3]) \to \Ext^2_\ce(K,K) = \ul{\ce}(K,S^2 K).
\]
Now since $\cc=\ul{\ce}$ is Calabi-Yau of CY-dimension $2$, we
have a canonical linear form
\[
\psi : \ul{\ce}(K,S^2 K) \to k
\]
and we define $\phi$ to be the restriction of $\psi$
to $\per(\cm)(P,P[3])$. Then
$\phi$ yields a morphism
\[
\per(\cm)(?, P[3]) \to D\per(\cm)(P,?)
\]
which sends $f$ to the map $g\mapsto \phi(fg)$. Clearly this is
a morphism of cohomological functors defined on the triangulated
category $\per(\cm)$. To check that it is an isomorphism, it suffices
to check that its evaluation at all shifts $S^i M^\wedge$, $i\in\Z$,
of representable functors $M^\wedge$ is an isomorphism. Now indeed, we
have
\[
\per(\cm)(S^i M^\wedge, P[3]) = \per(\cm)(S^i M^\wedge, Z[3]) =
\Ext^{3-i}(M^\wedge, Z).
\]
This vanishes if $i\neq 3$ and is canonically isomorphic to
\[
\Ext^1_\ce(M,K)
\]
for $i=3$. On the other hand, to compute $\per(\cm)(P, S^i M^\wedge)$, we have
to compute the homology of the complex
\[
0 \to \ce(M_0, M) \to \ce(M_1, M) \to \ce(M_2, M) \to \ce(M_3,M) \to 0.
\]
From the exact sequences
\[
0 \to \ce(M_0,M) \to \ce(M_1,M) \to \ce(K,M) \to \Ext^1_\ce(M_0,M)
\]
and
\[
0 \to \ce(K,M) \to \ce(M_2,M) \to \ce(M_3,M) \to \Ext^1_\ce(K,M)
\to \Ext^1_\ce(M_2,M) \ko
\]
using that $\Ext^1_\ce(M_0,M)=0=\Ext^1_\ce(M_2,M)$, we get that the
homology is $0$ except at $\ce(M_3,M)$, where it is
$\Ext^1_\ce(K,M_0)$.
Now we know that the canonical linear form on $\Ext^2_\ce(K,K)$ yields
a canonical isomorphism
\[
\Ext^1_\ce(M,K) \iso D\Ext^1_\ce(K,M)
\]
and one checks that it identifies with the given morphism
\[
\per(\cm)(M^\wedge[3], P[3]) \to D \per(\cm) (P, M^\wedge[3]).
\]
\end{proof}

The following corollary generalizes Proposition~6.2
of \cite{GeissLeclercSchroeer05a}. The fact that
$\mod\cm$ is abelian and of global dimension $3$
is a special case of the results of \cite{Iyama04a}.

\begin{corollary} Suppose that $\ce$ is abelian. Then
$\mod\cm$ is abelian of global dimension at most $3$ and for
each $X\in \mod\cm$ and each $Y\in\mod \ul{\cm}$, we have
canonical isomorphisms
\[
\Ext^i_{\mod\cm}(X,Y) \iso D\Ext^{3-i}_{\mod \cm}(Y,X) \ko i\in\Z.
\]
If $\ul{\cm}$ contains a non zero object, then $\mod\cm$
is of global dimension exactly $3$.
\end{corollary}

\begin{proof} Let $X$ be a finitely presented $\cm$-module
and let
\[
M_1^\wedge \to M_0^\wedge \to X \to 0
\]
be a presentation. We form the exact sequence
\[
0 \to K \to  M_1 \to M_0
\]
and then choose an exact sequence
\[
0 \to M_3 \to M_2 \to K \to 0
\]
using part a) of the proposition above. By splicing the two,
we obtain a complex
\[
0 \to M_3 \to M_2 \to M_1 \to M_0 \to 0 \ko
\]
whose image under the Yoneda functor is a projective resolution of
length at most $3$ of $X$. This shows that $\cm$
admits weak kernels
(hence $\mod \cm$ is abelian) and that $\mod\cm$ is of global
dimension at most $3$.
Thus the perfect derived category coincides with
the bounded derived category of $\mod\cm$. Now the claim about the
extension groups is obvious from the proposition.  For the last
assertion, we choose $X$ to be a non zero $\ul{\cm}$-module.
Then $\Ext^3_{\mod \cm}(X,X)$ is non zero.
\end{proof}

\section{Relative $(d+1)$-Calabi-Yau duality}
\label{s:relative-dp1}

\subsection{$d$-cluster tilting subcategories} \label{ss:setting}
Let $k$ be a field and $\cc$ a triangulated $k$-linear category with
split idempotents and suspension functor $S$. We suppose that all
$\Hom$-spaces of $\cc$ are finite-dimensional and that $\cc$ admits a
Serre functor $\Sigma$, \cf \cite{ReitenVandenBergh02}. Let $d\geq 1$
be an integer. We suppose that $\cc$ is Calabi-Yau of CY-dimension
$d$, \ie there is an isomorphism of triangle functors
\[
S^d \iso \Sigma.
\]
We fix such an isomorphism once and for all.

For $X,Y\in \cc$ and $n\in\Z$, we put
\[
\Ext^n(X,Y)=\cc\,(X,S^n Y).
\]
Assume that $\ct\subset\cc$ is a {\em $d$-cluster tilting subcategory}.
By this, we mean that $\ct$ is maximal $(d-1)$-orthogonal in the
sense of Iyama \cite{Iyama05}, \ie
\begin{itemize}
\item[a)] $\ct$ is a $k$-linear subcategory,
\item[b)] $\ct$ is functorially finite in $\cc$, \ie the
functors $\cc(?,X)|\ct$ and $\cc(X,?)|\ct$ are finitely generated
for all $X\in\cc$,
\item[c)] we have $\Ext^i(T,T')=0$ for all $T,T'\in\ct$ and
all $0< i <d$ and
\item[d)] if $X\in\cc$ satisfies $\Ext^i(T,X)=0$ for all
$0<i<d$ and all $T\in\ct$, then $T$ belongs to $\ct$.
\end{itemize}
Note that a), b), c) are self-dual and so is d) (by the
Calabi-Yau property). For $d=1$, we have $\ct=\cc$
by condition d). This definition is slightly different
from that in \cite{Thomas05} but most probably equivalent
in the context of [loc.cit.].
As in section~\ref{s:Gorenstein-property},
one proves that $\mod\ct$, the category of finitely presented
right $\ct$-modules, is abelian. Let
\[
F: \cc \to \mod\ct
\]
be the functor which sends $X$ to the restriction of $\cc(?,X)$
to $\ct$. For classes $\cu$, $\cv$ of objects of $\ct$, we denote
by $\cu*\cv$ the full subcategory of all objects $X$ of $\cc$
appearing in a triangle
\[
U \to X \to V \to SU.
\]
\begin{lemma} \label{lemma:presentation} Suppose that $d\geq 2$.
For each finitely presented module $M\in\mod\ct$, there
is a triangle
\[
T_0 \to T_1 \to X \to ST_0
\]
such that $FX$ is isomorphic to $M$. The functor $F$ induces an equivalence
\[
\cu/ (S\ct) \to \mod\ct \ko
\]
where $\cu=\ct*S\ct$ and $(S\ct)$ is the ideal of morphisms factoring
through objects $ST$, $T\in\ct$.
\end{lemma}

\begin{proof}
Since $\ct$ has split idempotents, the functor
$F$ induces an equivalence from $\ct$ to the category of projectives
of $\mod\ct$. Now let $P_1 \to P_0 \to M \to 0$ be a projective
presentation of $M$. Choose a morphism $T_1 \to T_0$ of $\ct$
whose image under $F$ is $P_1 \to P_0$. Define $X$ by the triangle
\[
T_1 \to T_0 \to X \to ST_1.
\]
Then, since $\Ext^1(T,T_1)$ vanishes for $T\in\ct$, we get an
exact sequence $FT_1 \to FT_0 \to FX \to 0$ and hence an isomorphism
between $FX$ and $M$. This also shows that $F$ is essentially
surjective since $X$ clearly belongs to $\cu$.
One shows that $F$ is full by lifting a morphism
between modules to a morphism between projective presentations
and then to a morphism between triangles. One shows that the
kernel of $F|\cu$ is $(S\ct)$ as in section~\ref{s:Gorenstein-property}.
\end{proof}

\subsection{Examples of $d$-cluster tilting subcategories}
One can use proposition~\ref{ss:d-tilting} below (\cf also
\cite{Wraalsen05} and \cite{Thomas05} \cite{FominReading05})
to produce the following examples: Consider the algebra $A$
given by the quiver
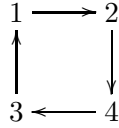
\[
\xymatrix{ 1 \ar[r] & 2 \ar[d] \\
3 \ar[u] & 4 \ar[l] }
\]
with the relations given by all paths of length $2$. Then $A$ is
$3$-cluster tilted of type $A_4$. It is selfinjective hence Gorenstein.
It is not hard to check that its
stable category is $4$-Calabi-Yau. More generally, a finite-dimensional
algebra of radical square $0$ whose quiver is an oriented cycle
with $d+1$ vertices is $d$-cluster-tilted of type $A_{d+1}$,
selfinjective and stably $(d+1)$-Calabi-Yau.

Let $A$ be the algebra given by the quiver
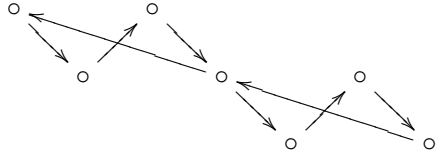
\[
\xymatrix@=0.5cm{
\bt \ar[dr] & & \bt \ar[dr] & &  & &  & \\
 & \bt \ar[ur] & & \bt \ar[dr] \ar[lllu]|-{}
& & \bt\ar[dr]  & & \\
 & &  & & \bt \ar[ru] & & \bt \ar[lllu] & }
\]
subject to the relations $\alpha\beta=0$ for all composable arrows
$\alpha$, $\beta$ which point in different directions
(the composition of the long skew arrows in non zero).
Then $A$ is $3$-cluster tilted of
type $A_{7}$. One can show that it is Gorenstein of dimension $1$
and that its stable Cohen-Macaulay category is $3$-Calabi-Yau.

\subsection{A counterexample} \label{ss:Iyamas-example}
The following example due to Osamu Iyama
shows that part d) of proposition~\ref{prop:gorenstein}
does not generalize from the $2$-dimensional to the $d$-dimensional case:
It is an example of a  $3$-cluster tilting subcategory of a $3$-Calabi Yau
triangulated category which is not Gorenstein.

Let $k$ be a field of characteristic $\neq 3$ and $\omega$ a primitive
third root of $1$. Let $S=k[[t,x,y,z]]$ and let the generator $g$ of
$G=\Z/3\Z$ act on $S$ by $gt=\omega t$, $gx=\omega x$, $gy=\omega^2y$
and $gz=\omega^2z$. Then the algebra $S^G$ is commutative,
an isolated singularity and Gorenstein of dimension $4$, 
its category of maximal Cohen-Macaulay modules
$\ce=\cmac(S^G)$ is Frobenius and the stable category $\cc=\ul{\cmac}(S^G)$
is $3$-Calabi-Yau. The $S^G$-module $S$ is maximal Cohen-Macaulay
and the subcategory $\ct$ formed by all direct factors of finite
direct sums of copies of $S$ in $\cc$ is a $3$-cluster tilting
subcategory. As an $S^G$-module, $S$ decomposes into the direct
sum of three indecomposable submodules $S_i$, $i=0,1,2$, formed
respectively by the $f\in S$ such that $gf=\omega^i f$. Note
that $S_0=S$ is projective. The full subcategory of $\cmac(S^G)$
with the objects $S_0$, $S_1$, $S_2$ is isomorphic to the
completed path category of the quiver
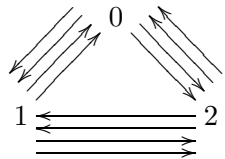
\[
\xymatrix{  & 0 \ar[rd] \ar@<1ex>[rd] \ar@<-2ex>[ld] \ar@<-3ex>[ld] &  \\
          1 \ar[ru] \ar@<1ex>[ru] \ar@<-2ex>[rr] \ar@<-3ex>[rr] &
        & 2 \ar[ll] \ar@<1ex>[ll] \ar@<-2ex>[lu] \ar@<-3ex>[lu] }
\]
subject to all commutativity relations, where, between any two vertices,
the arrows pointing counterclockwise are labelled $t$ and $x$ and the
arrows pointing clockwise $y$ and $z$. The vertex $0$ corresponds
to the indecomposable projective $S_0$, which vanishes in the
stable category $\cc$. Therefore the full subcategory $\ind(\ct)$ of
$\cc$ is given by the quiver
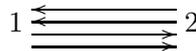
\[
\xymatrix{            1\ar@<-2ex>[rr] \ar@<-1ex>[rr] &
        & 2 \ar[ll] \ar@<-1ex>[ll]  }
\]
subject to the relations given by all paths of length $\geq 2$.
It is easy to see that the free modules $\ct(?,S_i)$, $i=1,2$, are
of infinite injective dimension. Thus, $\ct$ is not Gorenstein.

\subsection{Relative $(d+1)$-Calabi-Yau duality}
\label{ss:relative-dp1}
Let $d$ be an integer $\geq 1$.
Let $\ce$ be a $k$-linear Frobenius category with split idempotents.
Suppose that its stable category $\cc=\ul{\ce}$ has finite-dimensional
$\Hom$-spaces and is Calabi-Yau of CY-dimension $d$. Let $\ct\subset\cc$
be a $d$-cluster tilting subcategory and $\cm\subset\ce$ the preimage
of $\ct$ under the projection functor. In particular, $\cm$
contains the subcategory $\cp$ of the projective-injective
objects in $\cm$. Note that $\ct$ equals the quotient $\ul{\cm}$
of $\cm$ by the ideal of morphisms factoring through a projective-injective.

We call a complex {\em $\cm$-acyclic} if it is acyclic as a complex
over $\ce$ (\ie obtained by splicing conflations of $\ce$) and
its image under the
functor $X \mapsto \ce(?,X)|\cm$ is exact.

\begin{theorem}
\begin{itemize}
\item[a)] For each $X\in\ce$, there is an $\cm$-acyclic complex
\[
0 \to M_{d-1} \to M_{d-2} \to \ldots \to M_0 \to X \to 0
\]
with all $M_i$ in $\cm$.
\item[b)] For each finitely presented $\ul{\cm}$-module $Z$, there
is an $\cm$-acyclic complex
\[
0 \to M_{d+1} \to M_d \to \cdots \to M_1 \to M_0 \to 0
\]
together with isomorphisms
\[
Z \iso \cok(\ce(?,M_1) \to \ce(?,M_0)) \iso \Ext^1_{\ce}(?,Z_1) \ko
\]
where $Z_1$ is the kernel of $M_1 \to M_0$.

\item[c)] Let $Z$ be a finitely presented
$\ul{\cm}$-module. Then $Z$ considered as an $\cm$-module lies
in $\per(\cm)$ and we have a canonical isomorphism
\[
\per(\cm)(?,Z[d+1]) \iso D\per(\cm)(Z,?).
\]
\end{itemize}
\end{theorem}

\begin{proof} a) We choose a right $\ct$-approximation
$M_0' \to X$ in $\cc$ and complete it to a triangle
\[
Z'_0 \to M'_0 \to X \to SZ'_0.
\]
Now we lift the triangle to a conflation
\[
0 \to Z_0 \to M_0 \to X \to 0.
\]
By repeating this process we inductively construct conflations
\[
0 \to Z_i \to M_i \to Z_{i-1} \to 0
\]
for $1\leq i \leq d-1$. We splice these to obtain the complex
\[
0 \to Z_{d-1} \to M_{d-2} \to \cdots \to M_0 \to X \to 0.
\]
If follows from \ref{ss:resolutions} below that $M_{d-1}=Z_{d-1}$
belongs to $\cm$.

b) As in section~\ref{s:relative}, one sees that there is a
conflation
\[
0 \to Z_1 \to M_1 \to M_0 \to 0
\]
together with isomorphisms
\[
Z \iso \cok(\ce(?,M_1) \to \ce(?,M_0)) \iso \Ext^1_{\ce}(?,Z_1).
\]
Now we apply a) to $X=Z_1$.

c) Let $P$ be the complex constructed in b). We have
\[
\per(\cm)(P,P[d+1]) = \per(\cm)(P,Z[d+1]) =
\cok(Z(M_d) \to Z(M_{d+1})).
\]
Since $Z \iso \Ext^1_\ce(?,Z_1)$, we have to compute
\[
\cok (\Ext^1_\ce(M_d,Z_1) \to \Ext^1_\ce(M_{d+1},Z_1)).
\]
Put $X=SZ_1$.
We have $\Ext^1_\ce(?,Z_1) = \ul{\ce}(?,SZ_1)=\ul{\ce}(?,X)$ and
so we have to compute
\[
\cok( \ul{\ce}(M_d,X) \to \ul{\ce}(M_{d+1},X)).
\]
We now apply part~c) of proposition~\ref{prop:resolution} 
to the image in $\ul{\ce}$
of the complex
\[
0 \to M_{d+1} \to M_d \to \ldots \to M_d
\]
and the object $Y=Z_1=S^{-1}X$. In the notations of 
proposition~\ref{prop:resolution}, we obtain that
\[
\cok(GM_d \to GM_{d+1}) = GS^{-(d-1)}Y = \ul{\ce}(S^{-(d-1)}Y,?)
= \ul{\ce}(S^{-d}X,?).
\]
This yields in particular that
\[
\cok(\ul{\ce}(M_d,X) \to \ul{\ce}(M_{d+1},X)) =
\ul{\ce}(S^{-d}X, X) = \cc(X,S^d X).
\]
Thus, the canonical linear form on $\cc(X,S^d X)$ yields a morphism
\[
\per(\cm)(?,P[d+1]) \to D\per(\cm)(P,?).
\]
We have to check that it is an isomorphism on all representable
functors $M^\wedge[i]$, $i\in\Z$, $M\in\cm$. Now the group
\[
\per(\cm)(M^\wedge[i],P[d+1]) = \per(\cm)(M^\wedge, Z[d+1-i])
\]
vanishes if $i\neq d+1$ and equals $\ul{\ce}(M,X)$ for $i=d+1$.
On the other hand, we have to compute the group
\[
\per(\cm)(P,M^\wedge[i]).
\]
For $i=d+1$, it is isomorphic to $\ul{\ce}(S^{-d}X,M)$ by part c)
of proposition \ref{prop:resolution}. If we do not have
$2\leq i\leq d+1$, it clearly vanishes. For $2\leq i \leq d$, it
is the image of
\[
\ce(Z_{i-1}, M) \to \Ext^1(Z_{i-2},M).
\]
But $\Ext^1_\ce(Z_{i-1},M)=\ul{\ce}(Z_{i-2},SM)$ vanishes
because $Z_{i-2}$ is an iterate extension of objects in
\[
S^{-(i-2)}\ct, \ldots, S^{-1}\ct, \ct.
\]
\end{proof}

\subsection{Triangular resolutions} \label{ss:resolutions}
We work with the notations and assumptions of section~\ref{ss:setting}
and assume moreover that $d\geq 2$. Let $Y$ be an object of $\cc$.
Let $T_0 \to Y$ be a right $\ct$-approximation of $Y$. We define
an object $Z_0$ by the triangle
\[
Z_0 \to T_0 \to Y \to SZ_0.
\]
Now we choose a right $\ct$-approximation $T_1 \to Z_0$ and define
$Z_1$ by the triangle
\[
Z_1 \to T_1 \to Z_0 \to SZ_1.
\]
We continue inductively constructing triangles
\[
Z_i \to T_i \to Z_{i-1} \to SZ_i
\]
for $1<i\leq d-2$. By the proposition below, the
object $Z_{d-2}$ belongs to $\ct$. We put $T_{d-1}=Z_{d-2}$.
The triangles
\[
Z_i \to T_i \to Z_{i-1} \to S Z_i
\]
yield morphisms
\[
S^{-(d-1)}Y  \to S^{-(d-2)}Z_0 \to \cdots
\to S^{-1} Z_{d-3} \to Z_{d-2}
\]
so that we obtain a complex
\[
S^{-(d-1)}Y \to T_{d-1} \to T_{d-2} \to \ldots \to T_1 \to T_0 \to Y.
\]
Put $\ct\mbox{-}\mod=\mod(\ct\op)$.
Let
$
F: \cc \to \mod\ct \mbox{ and } G: \cc \to \ct\mbox{-}\mod
$
be the functors which take an object $X$ of $\cc$ to
$\cc(?,X)|\ct$ respectively $\cc(X,?)|\ct$.
The following proposition is related to Theorem~2.1 of
\cite{Iyama05}.

\begin{proposition} \label{prop:resolution} In the above notations, we have:
\begin{itemize}
\item[a)] The object $T_{d-1}=Z_{d-2}$ belongs to $\ct$.
\item[b)] The image of the complex
\[
T_{d-1} \to T_{d-2} \to \ldots \to T_1 \to T_0 \to 0
\]
under $F$ has homology $FY$ in degree $0$.
\item[c)] The image of the complex
\[
0 \to T_{d-1} \to T_{d-2} \to \ldots \to T_1 \to T_0
\]
under $G$ has homology $GS^{-(d-1)}Y$ in degree $d-1$.
If $Y$ belongs to $(S^{-1}\ct) * \ct$, then its homology
in degree $d-2$ vanishes.
\end{itemize}
\end{proposition}

\begin{proof} a) From the construction of the $Z_i$, we get that,
for $T\in\ct$, we have
\[
\Ext^1(T,Z_i)=0
\]
for all $i$ and that
\[
\Ext^j(T,Z_i) \iso \Ext^{j-1}(T,Z_{i-1}) \iso \cdots
\iso \Ext^1(T,Z_{i-j+1})=0
\]
for all $2\leq j \leq i+1$. Thus we have $\Ext^j(T,Z_{d-2})=0$
for $1\leq j\leq d-1$ and $Z_{d-1}$ indeed belongs to $\ct$.

b) This follows readily from the construction.

c) We consider the triangle
\[
S^{-1}T_{d-2} \to S^{-1}Z_{d-3} \to T_{d-1} \to T_{d-2}.
\]
For $T\in\ct$, the long exact sequence obtained by applying $\cc(?,T)$
combined with the vanishing of $\cc(S^{-1}T_{d-2},T)$ shows
that we have the isomorphism
\[
\cok(\cc(T_{d-2},T) \to \cc(T_{d-1},T)) \iso \cc (S^{-1}Z_{d-3},T).
\]
Now for $-1\leq i \leq d-4$, the maps $S^{-1}Z_i \to Z_{i+1}$ induce
isomorphisms
\[
\cc(S^{-(d-2)+i+1}Z_{i+1},T) \iso \cc(S^{-(d-2)+i}Z_{i},T)
\]
because of the triangles
\[
S^{i}T_{i+1} \to S^{i}Z_{i}\to  S^{i+1}Z_{i+1} \to
S^{i+1}T_{i+1}
\]
and the fact that
\[
\cc(S^{-(d-2)+i}T_{i+1},T)=0=\cc(S^{-(d-2)+i+1} T_{i+1},T).
\]
Here, we have taken $Z_{-1}=Y$. Finally, we obtain the isomorphism
\[
\cok(\cc(T_{d-2},T)\to \cc(T_{d-1},T))\iso \cc(S^{-(d-1)}Y,T).
\]

For the last assertion, let $T\in\ct$. By the sequence
\[
T_{d-1} \to T_{d-2} \to Z_{d-3} \to ST_{d-1}\ko
\]
the sequence
\[
\cc(Z_{d-3},T) \to \cc(T_{d-2},T) \to \cc(T_{d-1},T)
\]
is exact so that it suffices to show that
\[
\cc(T_{d-3}, T) \to \cc(Z_{d-3},T)
\]
is surjective. Thanks to the triangle
\[
S^{-1}Z_{d-4} \to Z_{d-3} \to T_{d-3} \to Z_{d-4} \ko
\]
it suffices to show that $\cc(S^{-1}Z_{d-4}, T)$ vanishes.
This is clear because the object $Z_{d-4}$ is an iterated
extension of objects in
\[
S^{-(d-2)}\ct, S^{-(d-3)}\ct, \ldots, S^{-1}\ct, \ct.
\]

\end{proof}

\subsection{$d$-cluster tilting categories from tilting subcategories}
\label{ss:d-tilting}
Let $\ch$ be an $\Ext$-finite hereditary $k$-linear Krull-Schmidt
category whose bounded derived category $\cd=\cd^b(\ch)$ admits a
Serre functor $\Sigma$.  Let $d\geq 2$ be an integer.
The {\em $d$-cluster category}
\[
\cc= \cd/(S^{-d}\circ \Sigma)^\Z \ko
\]
has been considered in \cite{Keller05} \cite{Wraalsen05} \cite{Thomas05}.
Let $\pi: \cd \to \cc$ be the projection functor. Let $\ct\subset\cd$ be
a tilting subcategory, \ie a full subcategory whose objects form a set
of generators for the triangulated category $\cd$ and such that
\[
\cd(T,S^i T')=0
\]
for all $T,T'\in\ct$ and all $i\neq 0$. It follows that there
is a triangle equivalence, \cf \cite{Keller91},
\[
\cd^b(\mod\ct) \to \cd^b(\ch)
\]
which takes the projective module $T^\wedge=\ct(?,T)$ in $\mod\ct$ to
$T$, for each $T\in\ct$. We assume that $\ct$ is locally bounded, \ie
for each indecomposable $T$, there are only finitely many
indecomposables $T'$ in $\ct$ such that $\ct(T,T')\neq 0$ or
$\ct(T',T)\neq 0$. We denote the {\em Nakayama functor} of $\mod\ct$
by $\nu$. Recall from section~3.7 of \cite{GabrielRoiter92} that it is
the unique functor endowed with isomorphisms
\[
D \Hom(P,M) = \Hom(M,\nu P)
\]
for all finitely generated projectives $P$ and all finitely presented
$\ct$-modules $M$. It follows that if we view a projective $P$
as an object of $\cd^b(\mod\ct)$, we have $\Sigma P \iso \nu P$.
The category of finitely generated projectives
is functorially finite in $\cd^b(\mod\ct)$: Indeed, if $M$ is
an object of $\cd^b(\mod\ct)$, we obtain a left approximation
by taking the morphism $P \to M$ induced by an epimorphism
$P \to H^0(M)$ with projective $P$, and we obtain a right
approximation by taking the morphism $M \to P'$ induced by
a monomorphism $H^0(\Sigma M) \to \nu P'$ into an injective
$\nu P'$. Therefore, $\ct$ is functorially finite in $\cd$.

The first part of the following proposition is proved in
section~3 of \cite{BuanMarshReinekeReitenTodorov04} for hereditary
categories with a tilting object for $d=2$, see also
Proposition~2.6 of \cite{Zhu05}. For arbitrary $d$, it
was proved in \cite{Wraalsen05} when $\ct$ is given by
a tilting module over a hereditary algebra. The second part has
been proved in \cite{AssemBruestleSchiffler06} for $d=2$ and
hereditary algebras.

\begin{proposition}
Assume that the objects of $\ct$ have their
homology concentrated in degrees $i$ with $-(d-2)\leq i \leq 0$.
Then $\pi(\ct)$ is a $d$-cluster tilting subcategory in $\cc$ and
for all objects $T,T'$ of $\ct$, we have a functorial isomorphism
compatible with compositions
\[
\pi(\ct)(\pi(T),\pi(T')) = \ct(T,T') \oplus \Ext^d_\ct(\nu T^\wedge,
T'^\wedge).
\]
\end{proposition}

Note that the proposition does not produce all $d$-cluster
tilting subcategories of $\cc$. For
example, as one checks easily, the non connected algebra $k\times k$ is
$3$-cluster tilted of type~$A_2$.

\begin{proof} Put $F=\Sigma^{-1} S^d$. For
$n\in\Z$, denote by $\cd_{\leq n}$ the full subcategory of
$\cd$ formed by the objects $X$ with $H^i(X)=0$ for $i>n$
and by $\cd_{\geq n}$ the full subcategory formed by the
$Y$ with $H^i(Y)=0$ for $i<n$. Then $\cd_{\leq 0}$ and $\cd_{\geq 0}$
are the aisles of the natural $t$-structure on $\cd$.
Thus the subcategory $\cd_{\geq 1}$ is the right perpendicular
subcategory of $\cd_{\leq 0}$ and we have in particular
$\Hom(X,Y)=0$ if $X \in\cd_{\leq 0}$ and $Y\in\cd_{\geq 1}$.
Since $\ch$ is hereditary, we also have $\Hom(X,Y)=0$ if $X\in\cd_{\geq 0}$
and $Y\in\cd_{\leq -2}$. We have $F \cd_{\leq 0} \subset \cd_{\leq -(d-1)}$.
Indeed, for $X\in\cd_{\leq 0}$ and $Y \in \cd_{\geq -(d-2)}$, we
have
\[
\cd(FX,Y)=\cd(\Sigma^{-1}S^d X,Y)=D\cd(Y,S^dX)=0
\]
because $S^dX$ belongs to $\cd_{\leq -d}$ and $Y$ to $\cd_{\geq
-d+2}$.

Let $X\in\cd$. Since we have $T\in\cd_{\leq 0} \cap \cd_{\geq
-(d-2)}$, there are only finitely many integers $n$ such that $\cd(T, F^n
X)\neq 0$ or $\cd(F^n X, T)\neq 0$ for some $T\in\ct$. Since $\ct$
is functorially finite in $\cd$, it follows that $\pi(\ct)$ is
functorially finite in $\cc$.

Let us show that for $T,T'\in\ct$, the space
$\cd(T,S^i F^n T')$ vanishes for $0<i$  and $n\geq 0$.
This is clear if $n=0$. For $n\geq 1$, we have
\[
S^i F^n T' \in \cd_{\leq -n(d-1)-i}
\mbox{ and } T \in \cd_{\geq -(d-2)}.
\]
Since we have
\[
-(d-2) + n(d-1)+i = (n-1)(d-1) +1 + i \geq 2
\]
it follows that $\cd(T,S^i F^n T')$ vanishes because $\ch$ is
hereditary.

Let us show that $\cd(F^n T, S^i T')$ vanishes for $i<d$ and $n>0$.
Let us first consider the case $n=1$. Then we have
\[
\cd(FT, S^i T') = \cd(\Sigma^{-1} S^d T, S^i T') =
D\cd(S^i T', S^d T) = D\cd(T', S^{d-i} T) =0
\]
since $i<d$. Now assume that $n>1$. Then we have
\[
F^n T \in \cd_{\leq -n(d-1)} \mbox{ and } S^i T' \in \cd_{\geq -(d-2) -i}.
\]
Since we have
\[
-(d-2)-i + n(d-1) = (n-1)(d-1) +1 -i \geq (d-1) + 1 - i = d-i >0 \ko
\]
it follows that $\Hom(F^n T, S^i T')$ vanishes.

We conclude that $\cd(T, S^i F^n T')$ vanishes for $0<i<d$ and
all $n\in\Z$ so that we have
\[
\cc(T, S^i T')=0
\]
for $0<i<d$. By re-examining the above computation, we also
find that $\cd(T, F^n T')$ vanishes for each integer $n\neq 0,1$ and that
we have
\[
\cd(T, FT') = \cd(T, \Sigma^{-1} S^d T') = \cd(\Sigma T, S^d T')
            = {\per(\ct)}(\Sigma T^\wedge, S^d T'^\wedge)
            = \Ext^d_\ct(\nu T^\wedge, T'^\wedge).
\]
Let $X\in\cd$. Let us show that if
\[
\Ext^i(\pi(T),\pi(X))=0
\]
for all $0<i<d$ and all $T\in\ct$, then $\pi(X)$ belongs to $\pi(\ct)$.
We may assume that $X$ is indecomposable and belongs to
\[
\cd_{\leq 0} \cap \cd_{\geq -(d-1)}.
\]
Let $\cu_{\leq 0}$, $\cu_{\geq 0}$ be the aisles obtained
as the images of the natural aisles in $\cd^b(\mod\ct)$ under
the triangle equivalence $\cd^b(\mod\ct)\iso \cd$ associated with $\ct$.
We claim that we have
\[
\cd_{\leq 0} \subset \cu_{\leq (d-1)}.
\]
Indeed, for $Y\in\cd_{\leq 0}$ and $i\geq d$, we have
$S^{-i}T \in \cd_{\geq 2}$ so that $\cd(S^{-i}T,Y)=0$
and this implies the claim.
Since $X$ lies in $\cd_{\leq 0}$ it follows that $X\in\cu_{\leq (d-1)}$.
Now the assumption on $X$ yields that
\[
\cd(S^{-i}T, X)=0
\]
for $0<i\leq d-1$. Thus $X$ lies in $\cu_{\leq 0}$.
Now we claim that we have
\[
\cd_{\geq 0} \subset \cu_{\geq 0}.
\]
Indeed, for $Y\in\cd_{\geq 0}$ and $i>0$, we have $\cd(S^i T, Y)=0$
since $S^i T \in\cd_{<0}$. This implies the claim.
Since $X$ lies in $\cd_{\geq -(d-1)}$, it follows that
$X$ lies in $\cu_{\geq -(d-1)}$.
Finally, for $T\in\ct$ and $0<i<d$, we have
\[
0 = \cd (T, S^i FX) = \cd(T, S^{i-d} \Sigma X) = D\cd(S^{i-d} X, T)
= D\cd(X, S^{d-i} T).
\]
Thus we have $\cd(X, S^i T)=0$ for $0<i<d$. But since $X$ lies
in $\cd_{\geq -(d-1)}$ and $T$ in $\cd_{\leq 0}$ we also have
$\cd(X, S^i T)=0$ for all $i\geq d$. Thus $X$ is left orthogonal
to $\cu_{\leq -1}$. Therefore, the object $Y \in \cd^b(\mod\ct)$
corresponding to $X$ via the triangle equivalence associated
with $\ct$ is in $\cd_{\leq 0}^b(\mod\ct)$ and left orthogonal to
$\cd_{\leq -1}^b(\mod\ct)$. Since $\mod\ct$ is of finite global
dimension, this implies that $Y$ is a projective
$\ct$-module. So $X$ lies in $\ct$.
\end{proof}


\def\cprime{$'$}
\providecommand{\bysame}{\leavevmode\hbox to3em{\hrulefill}\thinspace}
\providecommand{\MR}{\relax\ifhmode\unskip\space\fi MR }
\providecommand{\MRhref}[2]{%
  \href{http://www.ams.org/mathscinet-getitem?mr=#1}{#2}
}
\providecommand{\href}[2]{#2}

\end{document}